\documentclass{amsart}
\usepackage{amssymb,amsmath,mathrsfs,formel}
\usepackage[utf8]{inputenc}
\begin{document}
\def\sgn{\mathop{\mbox{{\rm sgn}}}}
\def\ggT{\mathop{\mbox{{\rm ggT}}}}
\def\Si{\mathop{\mbox{{\rm Si}}}}
\def\Re{\mathop{\mbox{{\rm Re}}}}
\def\cal{\mathscr}
\def\dd{\,{\mathrm d}}
\newtheorem{thm}{Theorem}
\newtheorem{lem}{Lemma}

\title[Rational points on the inner product cone]{Rational points on the inner product cone\\ via the hyperbola method}

\author{V. Blomer and J. Br\"udern}
\address{Mathematisches Institut, Bunsenstr. 3-5, 37073 G\"ottingen, Germany} \email{vblomer@math.uni-goettingen.de} \email{jbruede@mathematik.uni-goettingen.de}
\maketitle


\begin{abstract}
A strong quantitative form of Manin's conjecture is established
for a certain variety in biprojective space. The singular integral in an application of the circle method involves the third power of the integral sine function, and  is evaluated in closed form.
\end{abstract}

\section{Introduction.} 

\noindent
The natural inner product  on $3$-space defines an algebraic  variety 
\be{1} x_0y_0 + x_1 y_1 + x_2 y_2 =0. \ee
There are several ways to analyse  its diophantine properties quantitatively. The equation \rf{1} defines the isotropy cone of a senary quadratic form. In this view, one would ask for an asymptotic formula 
for the number $P(X)$ of solutions of \rf{1} in
integers not exceeding $X$ in modulus.  A prudent yet routine treatment of this classical question, by the circle 
method or otherwise, leads one to a positive number $C_0$ with
\be{00} P(X) = C_0 X^4 + O(X^{3}(\log X)^2). \ee

Alternatively, the left hand side of \rf{1} is a bilinear form in the variables  
${\bf x}=(x_0,x_1,x_2)$ and ${\bf y}=(y_0,y_1,y_2)$.  The first question then   is
to count the integral solutions of \rf{1} inside a box, with ${\bf x}$ and ${\bf y}$ bounded independently. Thus, writing
$|{\bf x}|=\max(|x_0|,|x_1|,|x_2|)$, an asymptotic expansion is desired for the number $M(X,Y)$ 
of ${\bf x},{\bf y}\in (\mathbb Z\backslash\{0\})^3$  that satisfy \rf{1} and lie inside the box
$|{\bf x}|\le X$, $|{\bf y}|\le Y$. The following theorem provides such a result  featuring the sums
\be{3.2}
A(n)=\sum_{j\le n}\frac{1}{j},\quad B(n)=\sum_{j\le n}\frac{1}{j^2},
\ee
and the function
\be{2}  F(n) = \Big(\frac{33}{2}-3B(n)\Big)n^2-\Big(\frac{21}{2}+3B(n)\Big)n+6A(n). 
\ee
Empty sums are to be read as $0$, and consequently, $F(0)=0$.

\begin{thm}\label{thm1}
 Let $\f32\le X\le Y$. Then
$$
M(X,Y)=4Y^2\sum^{\infty}_{q=1}\frac{\varphi(q)}{q}F\Big(\Big[\frac{X}{q}\Big]\Big)+O\big((XY)^{3/2} (\log X)\log Y\big).
$$
\end{thm}

Note that the reservoir of variables $\mathbf x, \mathbf y$ has about $(4XY)^3$ elements, so that the error term in Theorem \ref{thm1} corresponds to nearly square root cancellation. This seems to be the limit of what can be expected from
the circle method in this problem.  

The asymptotic expansion of $M(X,Y)$ is somewhat involved, but this is inevitable when $X$ is  small, and one insists on errors not much larger than $(XY)^{3/2}$. When $X$ gets larger,
the main term smoothes out. To see this, let  $t > 0$, write 
\begin{equation}\label{defg}
  G(t) = F([t]) - \frac{33 - \pi^2}{2} t^2,
\end{equation}
and note that
\begin{equation}\label{boundG}
  G(t) \ll \min(t, t^2).
\end{equation}
Then, equipped with the elementary formula 
 \be{4.7}
\sum^{\infty}_{q=1}\frac{\varphi(q)}{q^3}=\frac{\zeta(2)}{\zeta(3)}
\ee
we  readily find that
\begin{equation}\label{simple}
M(X,Y)=2(33-\pi^2)\frac{\zeta(2)}{\zeta(3)} (XY)^2  + O\left(XY^2 (\log X)\log Y\right)
\end{equation}
holds whenever $\f32 \leq X \leq Y$. This estimate exhibits square root cancellation in the shorter of the two variables, this time relative to the size of the leading term. We also note that $P(X)- M(X,X)$
is the number of solutions of \rf{1} with at least one variable $0$. Aided by 
Lemma 2 below, we easily find that $$P(X)= M(X,X)+O(X^3(\log X)^2), $$ 
and then  recover \rf{00} with $C_0= 2(33-\pi^2)\zeta(2){\zeta(3)}^{-1}.$

Continuing this line of thought, one may consider $x=(x_0:x_1:x_2)$ and $y=(y_0:y_1:y_2)$ as projective coordinates.
The equation \rf{1} then defines a variety $\cal V$ in biprojective space.  Note that
all $(x_0:x_1:x_2)\in\mathbb P^2(\mathbb Q)$ have exactly two representations by primitive ${\bf x}=(x_0,x_1,x_2)\in\mathbb Z^3$, and these differ only by sign, so that their norms $|{\bf x}|$ depend only on the projective point. If $(x,y)\in\mathbb P^2(\mathbb Q)\times P^2(\mathbb Q)$ is represented by  primitive ${\bf x}, {\bf y}\in\mathbb Z^3$ satisfying \rf{1}, one defines an anticanonical height of this point by
$
|{\bf x}|^2|{\bf y}|^2.
$
One is then interested in the number $N(B)$ of rational points on $\cal{V}$ of height not exceeding $B$. 
 Thus,
$4N(B)$ is the number of all primitive ${\bf x}, {\bf y}\in {\mathbb Z}^3$ satisfying \rf{1}, 
with
$|{\bf x}|^2|{\bf y}|^2\le B$. 

\begin{thm}\label{thm2}
 There is a real number $C$ such that one has
$$
N(B)=\frac{33-6\zeta(2)}{2\zeta(2)\zeta(3)} B\log B+ CB+O(B^{{7}/{8}}(\log B)^2).$$
\end{thm}

We shall prove Theorem 1 by the circle method, and deduce Theorem 2 from Theorem 1 by a development of Dirichlet's hyperbola method. 
A direct attempt to Theorem \ref{thm2}   by the circle method occurs in work of 
Spencer \cite{S}. He considered the Zariski open subset $\cal U$ of $\cal V$ where no projective coordinate vanishes. For the number $N_0(B)$ of all points in $\cal U$ that are counted by $N(B)$ he  obtained the asymptotic formula 
$
N_0(B)=c_0 B \log B +O(B),
$
with some positive $c_0$. It is somewhat unnatural to neglect the points on coordinate hyperplanes, and indeed, one has
\be{lem8} N(B) - N_0(B) =\Big(\f{48}{\zeta(3)} -\f{12}{\zeta(2)}\Big)B + O(B^{3/4}(\log B)^2). 
\ee
We shall confirm this in the final section of this paper. In particular, we see that Spencer's $c_0$ coincides with  
 the leading constant in Theorem \ref{thm2}. 

The height zeta function of $\cal V$ is defined by the series
\be{zeta}  \sum_{{\bf x}\cdot{\bf y}=0} (|{\bf x}||{\bf y}|)^{-2s}, \ee
with  summation restricted to primitive 
${\bf x}, {\bf y}\in{\mathbb Z}^3$. By partial summation and
Theorem 2, it follows that this series is absolutely convergent in $\Re s >1$ and admits an analytic continuation to the
region $\Re s > 7/8$, $s\neq 1$, with a double pole at $s=1$.  Moreover, for any $\theta>7/8$, in the region $\Re(s)\ge\theta$, $|s-1|\ge 1$, one  obtains the growth estimate $O(|s|)$  as a by-product. 

As an aside, we briefly comment on  an automorphic approach to the analysis of the variety $\cal{V}$. 
The solutions of \rf{1} in primitive vectors ${\bf x}, {\bf y}$ are in one-to-one correspondence with the cosets of $\Gamma_\infty\backslash {\rm SL}_3(\mathbb Z)$ where $\Gamma_\infty$ is the group of integral upper triangular unipotent matrices.  Therefore, a sum over the integral solutions of \rf{1} can be expressed in terms of the minimal parabolic Eisenstein series $E(g, s_1, s_2)$ on ${\rm GL}_3$, as defined in \cite[p.\ 100]{B}.  In particular (see \cite{B}, (7.2)), one finds that
$$ \sum_{{\bf x}\cdot {\bf y}=0} \big( (x_0^2+x_1^2+x_2^2)(y_0^2+y_1^2+y_2^2)\big)^{-s} = 
\textstyle E\big({\bf 1}, \f23 s, \f23 s\big) $$
where again the sum on the left runs over primitive ${\bf x}, {\bf y}$, and {\bf 1} denotes the 
$3\times 3$ identity matrix. The function  $ E({\bf 1}, \f23 s, \f23 s)$ is meromorphic in the complex plane, and in $\Re s >\f12$ the only  poles are at $s=1$ and $s=3/4$ (\cite{B}, Theorem 7.1). If one is prepared to replace the maximum norm
in the notion of height by the euclidean one, a conclusion similar to the one in Theorem 2 is readily recovered. One can express the height zeta function  
\rf{zeta} in terms of more general Eisenstein series associated to non-trivial ${\rm SO}_3$-types (cf.\ \cite{Bu}), but this leads to delicate convergence issues.  It should be noted  that the relation to the 
theory of Eisenstein series  plays a similar role, in much broader generality, in the important work of Franke, Manin and Tschinkel \cite{FMT}.

Despite its apparently more limited scope for the problem at hand, our more elementary treatment is of interest because a circle method approach promises success in situations where a suitable
group action  is not available. In a more general context of multi-homogeneous varieties, this was carried out in our  
forthcoming article \cite{BB}\footnote{There is a misprint in \protect\cite{BB}, (1.7).}. In the present paper, variations  of this circle of ideas are developed that provide more control on  the error term featured in Theorem \ref{thm2}. 
We refer to Section 5, perhaps the most original part of the present communication. Our approach differs  substantially from the version of the hyperbola method detailed in~\cite{BB}. 

As a by-product of our analysis, we are
 led to a curious identity involving the integral sine
$$
\Si  t=\int^{t}_0\frac{\sin\alpha}{\alpha}\dd \alpha
$$
that we have not been able to locate in the literature.

\begin{thm}\label{thm3} One has
\be{8}
\int^{\infty}_{0}\frac{(\Si t)^3}{t^3} \dd t=\frac{33}{32}\pi-\frac{1}{32}\pi^3.\ee
\end{thm}

\section{The circle method} 

The circle method argument prominently features the sum
\be{2.2}
f(\alpha)=\sum_{1\le |x|\le X}\sum_{1\le|y|\le Y}e(\alpha x y).
\ee

\begin{lem} Let $\frac{3}{2} \leq X \leq Y$. Let $\alpha\in\mathbb R$, and suppose that $a\in\mathbb Z$ and $q\in\mathbb N$ are coprime and satisfy $|q\alpha-a|\le q^{-1}$. Then
$$
f(\alpha)\ll (XYq^{-1}+X+q)\log Y.
$$
\end{lem}

{\it Proof.} If $q>XY$, then the estimate for $f(\alpha)$ supplied by Lemma 1 is weaker than the ``trivial'' bound $f(\alpha)\ll XY$ that is immediate from \rf{2.2}. Thus, we may suppose that $q\le XY$. By \rf{2.2},
$$
f(\alpha)\ll\sum_{1\le|x|\le X}\min(Y,\|\alpha x\|^{-1}),
$$
where as usual, $\|\beta\|$ denotes the distance of the real number $\beta$ to its nearest integer.
Lemma 2.1 of Vaughan \cite{hlm} now supplies an estimate for $f(\alpha)$ as is required. 

\begin{lem}
Let $\frac{3}{2} \leq X\leq Y$. Then
$$
\int^{1}_{0}|f(\alpha)|^2\dd \alpha\ll XY\log X.
$$
\end{lem}

{\it Proof.} By orthogonality, the integral in question equals the number of solutions of $xy=uv$ in integers $x,y,u,v$ satisfying
$$
 1\le |x|\le X,\quad 1\le |y|\le Y,\quad1\le |u|\le X, \quad 1\le |v|\le Y.
$$
By symmetry, it suffices to count solutions with $1\le u\le x\le X$. Write
$$
d=(x;u),\quad  x'=x/d,\quad u'=u/d.
$$
Then $(x';u')=1$, and hence $x'\mid v$, $u'\mid y$. Consequently, one has 
$v=x'w$, $y=u'w$ for some $w\in\mathbb Z$, and it follows that the integral in question does not exceed eight times the number of tuples $(d,x',u',w)$ with
$$
1\le d\le X,\quad  1\le u'\le x'\le X/d,\quad |w|\le Y/x'.
$$
Now $x' \leq X\le Y$, and  therefore,
$$
\int^{1}_{0}|f(\alpha)|^2\dd \alpha\ll\sum_{d\le X}\sum_{x\le X/d}x\sum_{w\le Y/x}1
\ll XY\log X.
$$
This completes the proof of Lemma 2.

\medskip

Recall that $M(X,Y)$ denotes the number of ${\bf x},{\bf y}\in (\mathbb Z\backslash\{0\})^3$ with $|{\bf x}|\le X$, $|{\bf y}|\le Y$ that satisfy \rf{1}. Then by \rf{2.2} and orthogonality,
\be{2.5}
M(X,Y)=\int^{1}_{0}f(\alpha)^3\dd \alpha. 
\ee
We  consider separately the contributions to \rf{2.5} that arise from major and minor arcs. To define the latter, let
\be{2.6}\textstyle
Q=\frac{1}{2}\sqrt{XY}.
\ee
The intervals $\{\alpha:|q\alpha-a|\le Q/(XY)\}$ with $1\le a\le q\le Q$ and $(a;q)=1$ are disjoint, and their union is denoted by $\frak M$. Let $\frak m=[Q/(XY),1+Q/(XY)]\backslash\frak M$.

By Dirichlet's theorem on diophantine approximation, for each $\alpha\in\mathbb R$ there are coprime integers $a$, $q$ with $1\le q\le 2\sqrt{XY}$ and $|q\alpha-a|\le Q/(XY)$. Hence $\alpha\in\frak m$ implies $q>Q$ here, and Lemma 1 yields
$$
\sup_{\alpha\in\frak m}|f(\alpha)|\ll XYQ^{-1}\log Y.
$$
By \rf{2.6} and Lemma 2, we now see that
\be{2.7}
\int_{\frak m}|f(\alpha)|^3 \dd \alpha\ll (XY)^{3/2}(\log X)\log Y. 
\ee

\medskip

We now develop an approximation to $f(\alpha)$ for use on the major arcs. When $\alpha\in\frak M$, there is a unique triple $(q,a,\beta)\in\NN^2\times \RR$ with 
$$ 1\le a\le q\le Q, \quad (a;q)=1,\quad  |\beta|\le Q/(qXY), \quad \al=\f{a}{q}+\beta.$$ 
We use this notation throughout the subsequent argument. 

We put $\al=(a/q)+\beta$ within \rf{2.2} and separate terms with $q\mid x$ from those with $q\nmid x$. In the portion where $q\mid x$, we restore the term $y=0$ to infer that
\be{2.8}
f(a/q+\beta)=f^*_q(\beta)+g_q(\alpha)+O(X)
\ee
where
$$
f^*_q(\beta)=\sum_{0<|x|\le X/q}\,\sum_{|y|\le Y} e(\beta qxy), \quad
g_q(\alpha)=\multsum{0<|x|\le X}{q\nmid x}\,\sum_{0<|y|\le Y} e(\alpha xy).
$$

\begin{lem}
 Let $\alpha$ be real, $a\in\mathbb Z$ and $q\in\mathbb N$ be coprime with $|q\alpha-a|\le Q/(XY)$. Then
$$
g_q(\alpha)\ll(X+q)\log q.
$$
\end{lem}

{\it Proof.} A very similar estimate occurs as Lemma 2.12 of Vaughan \cite{V2005}. We give the simple proof, for completeness. We have $g_1(\alpha)=0$ for all $\alpha$, whence we may suppose that $q\ge 2$. Summing over $y$ 
yields
$$
g_q(\alpha)\ll\multsum{0<|x|\le X}{q\nmid x}\|\alpha x\|^{-1}.
$$
Now $\|\alpha x\|\ge\|ax/q\|-|\beta x|$. But $q\nmid x$ and $(a;q)=1$ imply $\|ax/q\|\ge 1/q$, whereas $|\beta x|\le XQ/(qXY)\le 1/(2q)$ by \rf{2.6}.
This shows that $\|\alpha x\|\ge\frac{1}{2}\|ax/q\|$, and so,
$$
g_q(\alpha)\ll\multsum{0<|x|\le X}{q\nmid x}\Big\|\f{ax}{q}\Big\|^{-1}\ll \Big(\frac{X}{q}+1\Big)
\sum^{q-1}_{u=1}\Big\|\f{au}{q}\Big\|^{-1}.
$$
The conclusion of the lemma is now immediate.

\medskip

Recalling the relation between $\alpha\in\frak M$ and $a,q,\beta$, a function $f^*: {\frak M}\to \CC$ 
is defined via $f^*(\al)=f^*_q(\beta)$.  By \rf{2.8} and Lemma 3, we then have 
$$f(\alpha)=f^*(\al)+O(Q\log Q),$$
and hence, by binomial expansion,
\be{A1} 
f(\alpha)^3-f^*(\al)^3 \ll (|f(\alpha)|^2+|f^*(\al)|^2) Q\log Q.
\ee

\begin{lem} For $\frac{3}{2} \leq X \leq Y$, one has
$$ \int_{\frak M} |f^*(\al)|^2 \dd \al \ll XY \log X. $$
\end{lem}
Equipped with this lemma, we infer from  \rf{A1} and
Lemma 2 that 
\be{2.10}
\int_{\frak M}f(\alpha)^3\dd \alpha=
\int_{\frak M} f^{\ast}(\alpha)^3 \dd \alpha 
+O\big((XY)^{3/2}(\log X)\log Y\big).
\ee

The next tasks ahead of us are a proof of Lemma 4, and the evaluation of
the integral on the right hand side of \rf{2.10}. There are parallels in the treatment of the two problems. 
We begin with the obvious identity
$$
\sum_{|y|\le Y}e(\gamma y)= \frac{\sin\pi (2[Y]+1)\gamma}{\sin\pi\gamma},
$$
and apply it within \rf{2.8} to bring in the sum
\be{2.12}
w_q(\gamma)=\sum_{0<|x|\le X/q} \frac{\sin\pi(2[Y]+1)\gamma x}{\sin\pi\gamma x}
\ee
through the relation 
$ f^*(\al)= w_q(q\beta)$. 
Hence, on considering the contributions from the individual intervals comprising the major arcs separately, we find via the substitution $\gamma=q\beta$ that
\be{B18}  
 \int_{\frak M} |f^*(\al)|^2 \dd \al = \sum_{q\le Q}\frac{\varphi(q)}{q}\int^{Q/(XY)}_{-Q/(XY)}|w_q(\gamma)|^2 \dd\gamma\ee
and
\be{2.11}
\int_{\frak M}f^*(\alpha)^3\dd \alpha=\sum_{q\le Q}\frac{\varphi(q)}{q}\int^{Q/(XY)}_{-Q/(XY)}w_q(\gamma)^3\dd\gamma.
\ee

Next, we compare $w_q(\gamma)$ with the cognate sum
\be{2.13}
v_q(\gamma)=\sum_{0<|x|\le X/q} \frac{\sin\pi(2[Y]+1)\gamma x}{\pi\gamma x}.
\ee
With this in view, note that for $|t|\le\frac{\pi}{2}$, one has
$$
\frac{1}{\sin t}-\frac{1}{t}\ll |t|.
$$
Hence, for $|\gamma|\le (2X)^{-1}$, we may take $t=\pi\gamma x$ and sum over 
$0<|x|\le X/q$. By \rf{2.12} and \rf{2.13}, we infer that
\be{B19}
w_q(\gamma)=v_q(\gamma)+O(|\gamma|X^2q^{-2})
\ee
holds uniformly in $q\in\NN$. This implies the bound
$$ |w_q(\gamma)|^2 \ll |v_q(\gamma)|^2 + |\gamma|^2 X^4q^{-4}, $$
and from \rf{B18} and a routine estimation, we deduce that
\be{C18}
\int_{\frak M} |f^*(\al)|^2 \dd \al \ll \sum_{q\le Q}
\int^{Q/(XY)}_{-Q/(XY)}|v_q(\gamma)|^2 \dd\gamma + Q.
\ee

The trivial upper bounds
\be{2.14}
v_q(\gamma)\ll XY q^{-1}, \quad v_q(\gamma)\ll |\gamma|^{-1}\log X
\ee
are valid for $q\in\mathbb N$, $\gamma\in\mathbb R$, and are easily inferred from \rf{2.13}. It follows that $v_q\in L^2(\RR)$, whence $v_q$ has an $L^2$ Fourier transform $\hat v_q$, defined as the $L^2$ limit of the sequence of functions
$$ \int_{-n}^n v_q(\beta) e(-\al\beta)\dd\beta, $$
as $n\to\infty$. However, as is familiar, the limit
\be{B1} \lim_{n\to\infty}\int_{-n}^n \frac{\sin \pi t}{\pi t} e(-\al t) \dd t \ee
exists for all $\al\in\RR$, and equals the function $h:\RR \to [0,1]$
with $h(\al)=1$ for $|\al|< \f12$ and $h(\al)=0$ for $|\al|>\f12$. Hence, the
Fourier transforms of the individual summands in \rf{2.13} can be computed by rescaling \rf{B1}, and we find that
$$ \hat v_q (\al) = 2 \sum_{1\le x\le X/q} x^{-1} h\Big(\f{\al}{(2[Y]+1)x}\Big). $$
We apply Plancherel's theorem, asserting that
$$ \int_{-\infty}^\infty |v_q(\gamma)|^2 \dd\gamma = 
\int_{-\infty}^\infty |\hat v_q(\al)|^2 \dd \al, $$ to confirm the estimate
\begin{align} \int_{-\infty}^\infty |v_q(\gamma)|^2 \dd\gamma
&=  \sum_{1\le x,y\le X/q}\f{4}{xy}  \int_{-\infty}^\infty h\Big(\f{\al}{(2[Y]+1)x}\Big)h\Big(\f{\al}{(2[Y]+1)y}\Big) \dd \al\nonumber\\
& \le 8 (2[Y]+1)  \sum_{1\le x\le y\le X/q} \f1{y} \ll \f{XY}q.\label{vq}
\end{align}
Note that $v_q$ is an empty sum unless $q\le X$. Hence, we may sum the above over $q\le X$ and inject
 the result into \rf{C18}. This establishes Lemma 4.

\smallskip
For a similar treatment of \rf{2.11}, we write $w_q=v_q+(w_q-v_q)$ and then use
binomial expansion in conjunction with
\rf{B19} to see that
$$
w_q(\gamma)^3=v_q(\gamma)^3+O\big(|\gamma|^3X^6q^{-6}+|\gamma|X^2q^{-2}|v_q(\gamma)|^2\big)
$$
holds for $|\gamma|\le(2X)^{-1}$, and hence in particular when $|\gamma|\le Q/(XY)$. By \rf{vq}, integration over the latter range yields
$$
\int^{Q/(XY)}_{-Q/(XY)}w_q(\gamma)^3\dd\gamma =\int^{Q/(XY)}_{-Q/(XY)}v_q(\gamma)^3\dd\gamma +O\Big(
\f{Q^4X^2}{Y^{4}q^{6}}+ \f{QX^2}{q^3}\Big).
$$
 We now use the second part of \rf{2.14} together with \rf{vq} to control the error introduced by completing the integral
on the right to
\be{2.15}
J(q)=\int^{\infty}_{-\infty}v_q(\gamma)^3\dd\gamma, 
\ee
and then infer that
$$
\int^{Q/(XY)}_{-Q/(XY)}w_q(\gamma)^3\dd\gamma = J(q)
+O\Big( \f{X^4}{Y^{2}q^{6}}+ \f{QX^2}{q^3} + \f{X^2Y^2\log X}{qQ}\Big).
$$
We combine this with \rf{2.10} and \rf{2.11} to deduce that
\be{2.16}
\int_{\frak M}f(\alpha)^3\dd \alpha =\sum_{q\le Q}\frac{\varphi(q)}{q}J(q)+O\big((XY)^{3/2}(\log X)\log Y\big).
\ee
We wish to complete the sum over $q$ as well. By \rf{2.14}, \rf{vq} and \rf{2.15}, one finds that $J(q)\ll(XY)^2q^{-2}$, and hence,
$$
\sum_{q>Q}|J(q)|\ll(XY)^{2}Q^{-1} \ll (XY)^{3/2}.
$$
We may add this to \rf{2.16} and then add back in the contribution from the minor arcs. By \rf{2.7} and \rf{2.5}, we then conclude as follows.

\begin{lem}
 Let $\frac{3}{2}\le X \le Y$. Then
$$
M(X,Y)=\sum^{\infty}_{q=1}\frac{\varphi(q)}{q}J(q)+O\big((XY)^{3/2}(\log X)\log Y\big).
$$
\end{lem}

\section{ The singular integral} 
In this section, we compute the singular integral \rf{2.15}. Progress depends on the following identity that is reported in Gradshteyn-Ryzhik \cite{GR}, 3.763, no.\ 1 and 4. A straightforward proof is readily given via Cauchy's integral formulae.

\begin{lem}
 For a real number $\nu$, let $\Box (\nu)=\nu|\nu|$. Then, for positive real numbers $\omega_1,\omega_2,\omega_3$, one has
\begin{align*}
&\int\limits^{\infty}_{-\infty}\frac{\sin\omega_1t \,\sin\omega_2t\,\sin\omega_3t}{t^3}\dd t\\
=&\frac{\pi}{8}\left((\omega_1+\omega_2+\omega_3)^2+\Box(\omega_1-\omega_2-\omega_3)+\Box(\omega_2-\omega_3-\omega_1)+\Box(\omega_3-\omega_1-\omega_2)\right).
\end{align*}
\end{lem}

\medskip

We apply this result to \rf{2.15}, recalling \rf{2.13}. Then, using symmetry in $x_1,x_2,x_3$, one finds that
\be{3.1}
J(q)=(2[Y]+1)^2 S(X/q)
\ee
where 
\be{3.2B}
S(n)=\multsum{1\le x_i\le n}{i=1,2,3}\frac{(x_1+x_2+x_3)^2+3\Box(x_1-x_2-x_3)}{x_1x_2x_3}.
\ee
Recall the function $F$ as defined in \rf{2}.

\begin{lem}
 Let $n$ be a natural number. Then $S(n)=F(n)$.
\end{lem}

{\it Proof.}  In \rf{3.2B} we isolate terms where $\Box$ is non-negative. Then  
\be{3.2A}
S(n)=S_1(n)+6S_2(n)-3S_3(n)
\ee
where
\begin{align*}
S_1(n)&=\multsum{1\le x_i\le n}{i=1,2,3}\frac{(x_1+x_2+x_3)^2}{x_1x_2x_3},\\
S_2(n)&=\sum_{x_2+x_3\le x_1\le n}\frac{(x_1-x_2-x_3)^2}{x_1x_2x_3},\\
S_3(n)&=\multsum{1\le x_i\le n}{i=1,2,3}\frac{(x_1-x_2-x_3)^2}{x_1x_2x_3}.
\end{align*}

The goal is now to express these sums in terms of  $A(n)$ and $B(n)$, as defined in \rf{3.2}.
By symmetry in $x_1,x_2,x_3$, one readily confirms that
\be{3.3}
S_1(n)=\multsum{1\le x_i\le n}{i=1,2,3}\Big(\frac{3x_1}{x_2x_3}+\frac{6}{x_3}\Big)=\frac{3}{2}n(n+1)A(n)^2+6n^2A(n).
\ee
Similarly,
\be{3.4}
S_3(n)=\multsum{1\le x_i\le n}{i=1,2,3}\Big(\frac{3x_1}{x_2x_3}-\frac{2}{x_3}\Big)=\frac{3}{2}n(n+1)A(n)^2-2n^2A(n).
\ee
The sum $S_2$ is symmetric only in $x_2,x_3$. Hence, we only find that
\be{3.5}
S_2(n)=\sum_{x_2+x_3\le x_1\le n}\Big(\frac{x_1}{x_2x_3}+\frac{2x_3}{x_1x_2}-\frac{4}{x_3}+\frac{2}{x_1}\Big).
\ee

We treat the four summands separately. The rightmost term is
\be{3.6}
\sum_{x_2+x_3\le x_1\le n}\frac{2}{x_1}=\sum_{1\le x_1\le n}(x_1-1)=\frac{1}{2}n(n-1).
\ee
Similarly, the penultimate term contributes
\begin{align}
\sum_{x_2+x_3\le x_1\le n}\frac{4}{x_3}&=\sum_{x_3\le n}\frac{4}{x_3}\sum^{n}_{x_1=x_3}(x_1-x_3)\nonumber\\
&=\sum_{x_3\le n}\frac{2}{x_3}(n-x_3)(n-x_3+1)=2n(n+1)A(n)-3n^2-n. \label{3.7}
\end{align}

Next, one finds that
$$
\sum_{x_2+x_3\le x_1\le n}\frac{2x_3}{x_1x_2}=\sum_{x_2\le x_1\le n}\frac{(x_1-x_2+1)(x_1-x_2)}{x_1x_2}=T_1(n)+T_2(n)
$$
where
$$
T_j(n)=\sum_{x_2\le x_1\le n}\frac{(x_1-x_2)^j}{x_1x_2}.
$$
To evaluate $T_2(n)$, note that terms with $x_1=x_2$ do not contribute, and that the summands are invariant under $x_1\leftrightarrow x_2$. This shows that
$$
T_2(n)=\frac{1}{2}\multsum{x_1\le n}{x_2\le n}\frac{(x_1-x_2)^2}{x_1x_2}=\multsum{x_1\le n}{x_2\le n}\frac{x_1}{x_2}-n^2=\frac{1}{2}n(n+1)A(n)-n^2.
$$
One also has
\begin{align*}
T_1(n)&=\sum_{x_2\le x_1\le n} \frac{1}{x_2}-\sum_{x_2\le x_1\le n}\frac{1}{x_1}\nonumber\\
&=\sum_{x_2\le n}\frac{n-x_2+1}{x_2}-n=(n+1)A(n)-2n.\nonumber
\end{align*}
It follows that
\be{3.8}
\sum_{x_3x_2\le x_1\le n}\frac{2x_3}{x_1x_2}=\frac{1}{2}(n+1)(n+2)A(n)-n^2-2n.
\ee

For the first summand in \rf{3.5}, we first sum over $x_1$ to infer that
\begin{align*}
\sum_{x_3+x_2\le x_1\le n}\frac{x_1}{x_2x_3}&=\frac{1}{2}\sum_{x_2+x_3\le n}
\frac{n(n+1)-(x_2+x_3-1)(x_2+x_3)}{x_2x_3}\nonumber\\
&=\frac{1}{2}\big(n(n+1)U_0(n)+U_1(n)-U_2(n)\big),\nonumber
\end{align*}
where
$$
U_j(n)=\sum_{x+y\le n}\frac{(x+y)^j}{xy}.
$$
By symmetry,
$$
U_1(n)=\sum_{x+y\le n}\frac{2}{x}=2\sum_{x\le n}\frac{n-x}{x}=2(nA(n)-n)
$$
and
\begin{align*}
U_2(n)&=2\sum_{x+y\le n}\frac{x}{y}+n(n-1)=\sum_{y\le n}\frac{(n-y)(n-y+1)}{y}+n(n-1)\nonumber\\
&=n(n+1)A(n)-(2n+1)n+\frac{1}{2}n(n+1)+n(n-1)\nonumber\\
&=n(n+1)A(n)-\frac{1}{2}n^2-\frac{3}{2}n.\nonumber
\end{align*}
To compute $U_0(n)$, we substitute $z=x+y$ for $y$. This yields
$$
U_0(n)=\sum_{z\le n}\sum_{x\le z-1}\frac{1}{x(z-x)}
=\sum_{z\le n}\sum_{x\le z-1}\frac1z  \Big(\frac{1}{x}+\frac{1}{(z-x)}\Big), 
$$
and by symmetry, it follows that
$$
U_0(n)=\multsum{x\le n,z\le n}{x\ne z}\frac{1}{xz}=A(n)^2-B(n).
$$
The last four identities combine to
\be{3.9}
\sum_{x_3+x_2\le x_1\le n}\frac{x_1}{x_2x_3}=\frac{1}{2}n(n+1)\big(A(n)^2-A(n)-B(n)\big)+nA(n)+\frac{1}{4}n^2-\frac{1}{4}n.
\ee
By \rf{3.5}, \rf{3.6}, \rf{3.7}, \rf{3.8} and \rf{3.9},
\be{S2}
S_2(n)=\frac{1}{2}n(n+1)(A(n)^2-B(n))+(1-2n^2)A(n)+\frac{11}{4}n^2-\frac{7}{4}n.
\ee
The lemma now follows from \rf{3.2A}, \rf{3.3}, \rf{3.4} and \rf{S2}, on recalling \rf{2}.

\medskip

We are ready to deduce Theorem 1.  
 By \rf{3.1}, Lemma 7 and Lemma 5, one finds
$$
M(X,Y)=(2[Y]+1)^2\sum^{\infty}_{q=1}\frac{\varphi(q)}{q}F\Big(\Big[\frac{X}{q}\Big]\Big)+O\big((XY)^{3/2}(\log X) \log Y\big).
$$
Note that $(2[Y]+1)^2=4Y^2+O(Y)$. Since $F(n)\ll n^2$ holds for all $n$, we may indeed replace the first factor by $4Y^2$, with an error bounded by $O(X^2Y)$, and Theorem 1 follows.

\section{A curious identity} Theorem \ref{thm3}  is also readily available. In the work of Sections 2 and 3, we temporarily take $X=Y\in\NN$. Then, by
\rf{3.1}, Lemma 7, and \rf{defg} and \rf{boundG},
\be{5.1}
J(1)=(2X+1)^2F(X)=(66-2\pi^2)X^4+O(X^3).
\ee
Next, we derive an alternative expression for $J(1)$. Still restricted to the case $X=Y\in\NN$, we observe symmetry in the sum \rf{2.13} to infer that
$$ v_1(\gamma)=2\sum_{1\le x\le X}\frac{\sin\pi(2X+1)\gamma x}{\pi\gamma x}
=2\sum_{0\le x\le X}\frac{\sin\pi(2X+1)\gamma x}{\pi\gamma x} + O(X), $$
with the obvious interpretation that the function $(\sin t)/t$ takes the value $1$ at $t=0$. 
We now apply Euler's summation formula. When $|\gamma|\le X^{-3/2}$, this gives
$$
v_1(\gamma)=2\int^X_0\frac{\sin\pi(2X+1)\gamma t}{\pi\gamma t}\dd t+O(X^{3/2})=\frac{2}{\pi\gamma}\Si\,(\pi(2X+1)X\gamma) +O(X^{3/2}).
$$
By \rf{2.14} and \rf{2.15}, we have
$$ J(1) = \int^{X^{-3/2}}_{-X^{-3/2}}v_1(\gamma)^3\dd\gamma+O(X^{7/2}), $$
and we then insert the expansion of $v_1(\gamma)$ to deduce that
\begin{align}
J(1)&=\frac{8}{\pi^3}\int^{X^{-3/2}}_{-X^{-3/2}}\big(\Si\,( \pi(2X+1)X\gamma)\big)^3\,\frac{\mathrm d\gamma}{\gamma^3}+O(X^{7/2})\nonumber\\
&=\frac{32}{\pi}X^4\int^{\infty}_{-\infty}\frac{(\Si\gamma)^3}{\gamma^3}\dd\gamma+O(X^{7/2}).\label{5.2}
\end{align}
Since  $(\Si\gamma)^3\gamma^{-3}$ is even,
the identity \eqref{8} follows from comparing \rf{5.1} with \rf{5.2}.

\section{ The hyperbola method} The transition from an asymptotic formula for $M(X,Y)$ to one for $N(B)$ is the theme of this and the following section. Let $M'(B)$ denote the number of ${\bf x},{\bf y}\in (\mathbb Z\backslash\{0\})^3$ satisfying $(1)$ and $|{\bf x}|^2|{\bf y}|^2\le B.$  An asymptotic formula for $M'(B)$ is given in \eqref{4.11} below.
 One may  apply Dirichlet's hyperbola method in its classical form to the natural numbers $|{\bf x}|$ and $|{\bf y}|$. This yields
$$
M'(B)=M(B^{1/4}, B^{1/4})+2\sum^{[B^{1/4}]}_{l=1}\left(M\Big(l,\f{B^{1/2}}{l})-M\Big(l,\f{B^{1/2}}{l+1}\Big)\right).
$$
Unfortunately, in this form the method is not particularly useful, because each of the differenced terms within the sum over $l$ will import an error of size $B^{3/4}(\log B)^2$ from Theorem 1,  summing up to $B(\log B)^2$, far too big for a successful estimation. To surmount this difficulty, one works with fewer terms. In \cite{BB} we proposed a combination of linear and geometric progressions as sample points, while here we employ a quadratic sequence. Define $L \in \Bbb{N}$ by $(L-1)^2 < B^{1/4} \leq L^2$, so that
\begin{equation}\label{BL}
B^{1/8} \ll L \ll  B^{1/8} \quad \text{and} \quad L^2 = B^{1/4} + O(B^{1/8}).
\end{equation}
 Then, one readily confirms the lower bound
\begin{equation}\label{4.3}
M'(B)\ge M(B^{1/2}L^{-2},B^{1/2}L^{-2})+2\Xi \end{equation}
where 
$$ \Xi=\sum^{L-1}_{l=1}\left(M\Big(l^2,\f{B^{1/2}}{l^2}\Big)-M\Big(l^2,\f{B^{1/2}}{(l+1)^2}\Big)\right).
$$
Likewise, one finds
 the corresponding upper bound
\begin{equation}\label{4.4}
M'(B)\le M(B^{1/4},B^{1/4})+2\sum^L_{l=2}\left(M\Big(l^2,\f{B^{1/2}}{(l-1)^2}\Big)
-M\Big(l^2,\f{B^{1/2}}{l^2}\Big)\right).
\end{equation}

We now evaluate the right hand sides of \rf{4.3} and \rf{4.4} with the aid of Theorem~\ref{thm1}. It will be convenient to write
\be{4.5}
c=66-2\pi^2,
\ee
so that we conclude from \eqref{simple} and \eqref{BL} that
$$
M(B^{1/2}L^{-2},B^{1/2}L^{-2})=\frac{\zeta(2)}{\zeta(3)}cB+O (B^{7/8}  ),
$$
and the same asymptotic formula holds for $M(B^{1/4}, B^{1/4})$. 
Similarly, one has
\be{4.8}
\Xi =4B\sum^{L-1}_{l=1}\sum^\infty_{q=1}\frac{\varphi(q)}{q}F\Big(\Big[\frac{l^2}{q}\Big]\Big)\Big(\frac{1}{l^4}-\frac{1}{(l+1)^4}\Big)+O\big(B^{7/8}(\log B)^2\big).
\ee
By \eqref{defg}, the main term on the right hand side   may now be written in the form
\be{extra}
cB\sum_{l=1}^{L-1}\sum^\infty_{q=1}\frac{\varphi(q)}{q^3}\Big(1 - \frac{l^4}{(l+1)^4}\Big)+B\sum_{l=1}^{L-1}\sum^\infty_{q=1}\frac{\varphi(q)}{q}G\Big(\frac{l^2}{q}\Big)\Big(\frac{1}{l^4}-\frac{1}{(l+1)^4}\Big).
\ee
In the first summand here, we apply \rf{4.7} and then use Euler's summation formula to find a real number $c_1$ with   
$$\sum_{l=1}^{L-1}\Big(1 - \frac{l^4}{(l+1)^4}\Big) = 4 \log L + c_1 + O(L^{-1}) = 2 \log B^{1/4} + c_1 + O(B^{-1/8}).$$
Hence the first summand in \rf{extra} equals
  $2\zeta(2)\zeta(3)^{-1}cB\log B^{1/4}+c'B+O(B^{7/8})$, for some suitable $c'\in\mathbb R$. 
    In the second summand in \rf{extra}, we sum over all $l$ to remove dependence on $L$. By \rf{boundG}, this introduces an error not exceeding
$$
B\sum_{l>L}\sum^\infty_{q=1}\min\Big(\frac{1}{ql^3},\frac{1}{lq^2}\Big)\ll BL^{-2}\log L \ll B^{3/4} \log B.
$$
Collecting together, we now deduce from \rf{4.8} that
$$
\Xi=2\frac{\zeta(2)}{\zeta(3)}cB\log B^{1/4}+c''B+O\big(B^{7/8}(\log B)^2\big).
$$
holds with some suitable real number $c''$, and then we obtain the lower bound
\be{4.10}
M'(B)\ge 4\frac{\zeta(2)}{\zeta(3)}cB\log B^{1/4} +c'''B
+O\big(B^{7/8}(\log B)^2\big),
\ee
with some $c'''\in\mathbb R$.  One may apply the same argument to the right hand side of \rf{4.4}, and with hardly any change in the preceding computation, this leads to an upper  bound for $M'(B)$ in which the leading term coincides with the one in \rf{4.10}. Thus, we have established the asymptotic formula
\be{4.11}
M'(B)=\frac{\zeta(2)}{\zeta(3)}cB\log B+c'''B+O\big(B^{7/8}(\log B)^{2}\big).
\ee

\section{Proof of Theorem \ref{thm2}}

Recall that $N_0(B)$ is the number of rational points of height not exceeding $B$  on the Zariski open subset $\cal U$ of $\cal V$ where no projective coordinate vanishes. According to a comment preceding the statement of Theorem 2, one observes that $4N_0(B)$ is the number of {\it primitive} vectors ${\bf x}, {\bf y}\in\mathbb Z^3$ counted by $M'(B)$. We remove the conditions $(x_0;x_1;x_2)=(y_0;y_1;y_2)=1$ by Möbius inversion to arrive at
$$
4N_0(B)=\sum_{nm\le B^{1/2}}\mu(n)\mu(m)M'(B/(nm)^2).
$$
By \rf{4.11} and straightforward estimates, it follows that there is a number $C'$ with
$$
N_0(B)=\frac{c}{4\zeta(2)\zeta(3)}B\log B+C'B+O\big(B^{7/8}(\log B)^{2}\big).
$$
By \rf{4.5} and \rf{lem8}, the conclusion of Theorem \ref{thm2} is now available.

\medskip

It remains to establish \rf{lem8}. As above, we see that $4(N(B) - N_0(B))$ equals 
 the number of  primitive vectors $\mathbf x, \mathbf y \in \Bbb{Z}^3$ with $|{\bf x} |^2 | {\bf y}|^2 \leq B$ satisfying \eqref{1} and $x_0x_1x_2y_0y_1y_2 = 0$. Let $W_j$ denote the number of pairs $\mathbf x, \mathbf y$ counted here, with exactly $j$ of the integers $x_0,x_1,x_2,y_0,y_1,y_2$ zero. Since a primitive $\mathbf x\in\ZZ^3$ is non-zero, at most two of $x_0,x_1,x_2$ are zero. Further, if (say) $x_0=x_1=0$ then $\mathbf x$ is primitive if and only if $x_2=\pm 1$. Thus  $W_4\ll 1$ and $W_j=0$ for $j\ge 5$, and hence,
\be{002} 4(N(B) - N_0(B)) = W_1+W_2+W_3 +O(1). \ee

For the rest of this section we put $Z=B^{1/2}$ in the interest of readability.
For a solution counted by $W_3$, one of $\bf x, \bf y$ must have two of its coordinates zero. Hence, we may suppose temporarily that $x_0=x_1=0$. Then, as just observed, $x_2=\pm 1$, and \rf{1} implies $y_2=0$. Because there are two choices for $x_2$ and six symmetric choices for the particular role attributed to $x_0,x_1,y_2$, it follows that
$$ W_3 = 12 \multsum{0<|y_0|, |y_1|\le Z }{(y_0;y_1)=1} 1 = \f{48}{\zeta(2)}Z^2 + O(Z \log Z). $$
Here we have used Satz 1.3.4 and the discussion following that result in Br\"udern~\cite{JB}.

Next, consider a solution counted by $W_2$. Suppose that $x_0=0$, say. If it were the case that $x_1y_1=0$ then by \rf{1} one would have $x_2y_2=0$, and such a solution is not counted by $W_2$. We conclude that $x_0=y_0=0$, and \rf{1} reduces to $x_1y_1+x_2y_2=0$. But $\bf x,y$ are primitive, and we conclude that $(x_1,x_2)=\pm (y_2, -y_1)$. The size constraint
$|{\bf x} | | {\bf y}| \leq Z$ implies $x_1^2 \le Z$ and $x_2^2 \le Z$, and observing symmetry, it follows that
 $W_2 \ll Z$.

This leaves us with the estimation of $W_1$. By symmetry, $W_1= 6 W'_1$ where $W'_1$ is the number of solutions counted by $W_1$ with $x_0=0$. Hence $W'_1$ is the number of integers $x_1,x_2,y_0,y_1,y_2$, all non-zero,   with $x_1y_1+x_2y_2=0$ and
$$ (x_1;x_2)=(y_0;y_1;y_2)=1,\quad |x_iy_j|\le Z \quad (1\le i\le 2,\, 0\le j\le 2). $$
For any tuple counted here, we have $x_1\mid y_2$ and $x_2\mid y_1$. Hence we may write $y_2=ux_1$ and $y_1 = -ux_2$, with some non-zero integer $u$, and it follows that $W'_1$ is the number of non-zero integers $x_1,x_2,u$ and $y=y_0$ with
$$   (x_1;x_2)=(y;u)=1, \quad |x_iy| \le Z, \quad |u|x_i^2 \le Z \quad (i=1,2). $$ 
Thus, $W'_1 = 16 W''_1$ where $W''_1$ is number of positive $x_1,x_2,u,y$ counted by
$W'_1$. 

Note that $(x_1;x_2)=1$ and $x_1=x_2>0$ imply that $x_1=x_2=1$. Let $W_1^+$ denote the number of tuples counted by $W''_1$ where $x_1<x_2$, and let $W_1^0$ denote the number of such tuples with $x_1=x_2=1$. Then, by symmetry in $x_1$ and $x_2$, we have $W''_1 = W_1^0 + 2W_1^+$. Moreover, as in the treatment of $W_3$, we see that
$$W_1^0= \multsum{1\le u,y\le Z}{(u;y)=1} 1 = \f{1}{\zeta(2)}Z^2 + O(Z \log Z). $$
For the evaluation of $W_1^+$, note that 
for a given $x=x_2$, there are $\phi(x)$ choices for $x_1$. Further, the conditions on the natural number $y$ read
$y\le Z/x$ and $(y;u)=1$. As is well known, there are $(Z/x)(\phi(u)/u) + O(d(u))$ such $y$; here
 $d(u)$ is the number of divisors of $u$. We now see that $W_1^+$ equals
$$  \multsum{ux^2\le Z}{x\ge 2} \phi(x)\Big(\f{Z\phi(u)}{xu} + O(d(u)) \Big) =
Z \sum_{2\le x \le \sqrt Z} \f{\phi(x)}{x} \sum_{u\le Z/x^2} \f{\phi(u)}{u} + O(Z(\log Z)^2). $$
The remaining sums are routinely computed via  \cite[Satz 1.3.4]{JB} and partial summation,
and one then has 
$$  W_1^+=\f{Z^2}{\zeta(2)}\sum_{x=2}^\infty \f{\phi(x)}{x^3} +O\big(Z^{3/2}(\log Z)^2\big)
= Z^2 \Big(\f{1}{\zeta(3)} - \f{1}{\zeta(2)}\Big) +O\big(Z^{3/2}(\log Z)^2\big). $$
Collecting together, we first find that
$$ W_1 = 96\Big(\f{2}{\zeta(3)} -\f{1}{\zeta(2)}\Big) Z^2 +  O\big(Z^{3/2}(\log Z)^2\big), $$
and then insert our findings on $W_j$ in \rf{002} to confirm \rf{lem8}.

\end{document}